\documentclass[11pt]{article}
\usepackage{amssymb,amsfonts,amsmath,latexsym,epsf,tikz,url}
\newtheorem{theorem}{Theorem}[section]

\newtheorem{conjecture}[theorem]{Conjecture}
\newtheorem{corollary}[theorem]{Corollary}
\newtheorem{lemma}[theorem]{Lemma}
\newtheorem{rem}[theorem]{Remark}
\newtheorem{example}[theorem]{Example}
\newtheorem{definition}[theorem]{Definition}

\newcommand{\proof}{\noindent{\bf Proof.\ }}
\newcommand{\qed}{\hfill $\square$\medskip}

\begin{document}

\title{A graph which recognizes idempotents of a commutative ring}

\author{Hamid Reza Dorbidi$^{a}$ and
Saeid Alikhani$^{b,}$\footnote{Corresponding author}  }

\date{\today}

\maketitle

\begin{center}

    $^a$Department of Mathematics, Faculty of Science, University of Jiroft, P.O. Box 78671-61167, Jiroft, Iran\\
    \medskip

$^b$Department of Mathematics, Yazd University, 89195-741, Yazd, Iran\\
{\tt hr\_dorbidi@ujiroft.ac.ir, alikhani@yazd.ac.ir}\\

\end{center}


\begin{abstract}
In this paper we introduce and study a graph  on the set of ideals of a
commutative ring $R$. The vertices of this graph are non-trivial ideals of $R$ and
two distinct ideals $I$ and $J$ are adjacent if  and only  $IJ=I\cap J$. We
obtain some properties of this graph and study its
relation to the structure of $R$.
\end{abstract}

\noindent{\bf Keywords:} Graph, Diameter, Ring, Idempotent

\medskip
\noindent{\bf AMS Subj.\ Class.:} 05C25,20F65

\section{Introduction}

The study of algebraic structures, using the properties of graph theory, tends to exciting research topic in the last years. There are many papers on assigning a graph to a ring. Also some graph structures on the set of ideals of a ring $R$  are defined in the last decade.
The intersection graph of a ring (\cite{chak})  is a graph whose its vertices are non trivial ideals of $R$ and two distinct vertices $I$ and $J$ are adjacent if and only if $I\cap J\neq 0$. This graph is denoted by $\Gamma(R)$.
Authors in \cite{chak}, have characterized the rings $R$ for which the graph $\Gamma(R)$ is connected and obtained several necessary
and sufficient conditions on a ring $R$ such that $\Gamma(R)$ is a complete graph.  Also they determined the values of $n$ for which the graph of $\mathbb{Z}_n$ is Eulerian and Hamiltonian. Akbari, et al. in  \cite{akb}   determined all rings whose clique number of the intersection graphs of ideals is finite. Also  they showed that, if the clique number of $\Gamma(R)$ is finite, then its  chromatic number is finite and if $R$ is a reduced ring, then both are equal.

 Annihilating ideal graph (\cite{beh}) is a graph which its  vertices  are ideals with nonzero annihilators and two distinct vertices $I$ and $J$ are adjacent if and only if
$IJ=0$. This graph is denoted by $AG(R)$. It is shown that (\cite{beh}) if $R$ is not a domain, then $AG(R)$  has ascending chain condition (respectively, descending chain condition) on vertices if and only if $R$ is Noetherian (respectively, Artinian). Also the connectivity, the diameter and  coloring  of $AG(R)$ has studied in \cite{beh}.

Comaximal ideal graph of a ring is defined in \cite{wu}. The vertices  are ideals which are not contained in the
Jacobson  of $R$ and two distinct vertices $I$ and $J$ are adjacent if and only if
$I+J=R$. Also this graph has studied in \cite{hrdor}.

 Intersection graph is the complement of a zero divisor graph of a semigroup and annihilating ideal graph and
comaximal ideal graph are
the zero divisor graph of some semigroups. So these graphs share many properties with zero divisor graphs. In this paper we study a new graph on the set of non-trivial ideals of a commutative
ring. Also we study the relationship between the primary decomposition of ideals of a ring  $R$ and connectivity of new graph.
\medskip

First we would like to recall some facts and notations related to this paper.
Let $G=(V,E)$ be a graph. The order of $G$ is the number
of vertices of $G$. For two graphs $G_1 = (V_1, E_1)$ and $G_2 = (V_2, E_2)$, the disjoint union of $G_1$ and $G_2$ denoted by $G_1\cup G_2$ is the graph
with vertex set $V_1\cup V_2$ and edge set $E_1 \cup E_2$.
 A graph without edges is called an empty (null) graph. If every two distinct vertices
of a graph of order $n$ are adjacent, the graph is called a complete graph and is denoted by $K_n$. A clique of a graph $G$ is a complete subgraph of $G$ and clique number of $G$  is the number of vertices in a maximum clique of $G$.
For every vertex $v\in V$,
the open neighborhood of $v$ is the set $N(v) =\{u\in V: uv \in E\}$ and the closed neighborhood is the set $N[v] = N(v)\cup \{v\}$. For every vertex $v\in V(G)$,
the degree of $v$ is $|N(v)|$, i.e., the number of edges incident with $v$.
 Let $G_1,\ldots,G_k$ be some graphs. Then we consider  $G_1\times\cdots\times G_k$ as  a graph whose vertex set is $V(G_1)\times\cdots\times V(G_k)$ and two vertices $(v_1,\ldots,v_k)$ and $(u_1,\ldots,u_k)$ are adjacent if and only if $v_i$ and $u_i$ are adjacent in $G_i$ for each $i$. As usual we show the distance between two vertices $v$ and $w$, by
 $d(v,w)$. The eccentricity $\epsilon(v)$ of a vertex $v$ is the greatest  distance between $v$ and any other vertex.
The diameter of a graph $G$ is denoted by $diam(G)$ and  is the maximum eccentricity of any vertex in the graph.

\medskip

Throughout this paper all rings are commutative with unit
element. A Von Neumann regular ring is a ring $R$ such that for
every $a\in R$ there exists an $x\in R$ such that $a=axa$. This
implies that $(ax)^2=ax$ and $\langle a\rangle=\langle
ax\rangle$. So every principal ideal is generated by an
idempotent element. If for any nonzero ideal $J$ of $R$, we have
$I\cap J\neq 0$ we say that $I$  is a large ideal.  A ring with a
unique maximal ideal is called a local ring. We denote the set of
all maximal ideals and all prime ideals of $R$ by $Max(R)$ and
$Spec(R)$, respectively. Also $M(R)$ denotes the set of minimal ideals of $R$. The intersection of all maximal ideals
of $R$ is called the Jacobson radical of $R$ and is denoted by
$J(R)$. The intersection of all prime ideals is the set of all
nilpotent elements and is denoted by $Nil(R)$. It is clear that
$Nil(R)\subseteq J(R)$.  A discrete valuation ring (DVR) is a principal ideal domain (PID) with exactly one non-zero maximal ideal.

The radical of an ideal $I$ is denoted by $r(I)$ and is defined as $\{r\in R:a^n\in I\}$. It is a standard fact that $r(I)=\bigcap_{I\subseteq P} P$.
An ideal $Q\neq R$ is called a primary ideal if $ab\in Q$ implies $a\in Q$ or $b\in r(Q)$.
It is easily seen that if $Q$ is a primary ideal then $r(Q)$ is a prime ideal.
Also if $r(Q)$ is a maximal ideal then $Q$ is a primary ideal.
We say that an ideal $I$ has a primary decomposition if $I=\bigcap_{i=1}^n Q_i$ where $Q_i$ are primary ideals. The set $\{r(Q_i)\}$ is denoted by $Ass(I)$ and is called the set of associated prime ideals of $I$. For two ideals $I$ and $J$ of $R$, we denote the set $\{r\in R: rJ\subseteq I\}$ by $(I:J)$.
Also, we denote the finite field with $q$ elements by $\mathbb{F}_q$.

\section{Introduction to a new graph }

In this section, we introduce a new graph on the set of ideals of a commutative ring $R$ which we denote it by $\Gamma_0(R)$ and study its properties.

\begin{definition}
Let $R$ be a commutative ring.  The vertices of the  graph $\Gamma_0(R)$  are non-trivial ideals of $R$ and
two distinct ideals $I$ and $J$ are adjacent if and only if  $IJ=I\cap J$.
\end{definition}
We also need the following definition in some cases:
\begin{definition}
Let $R$ be a commutative ring.  The vertices of the  graph $\Gamma_1(R)$ are all  ideals of $R$ and
two  ideals $I$ and $J$ are adjacent if and only if  $IJ=I\cap J$.
\end{definition}

By this definition, there is a loop in the  vertex $I$ of $\Gamma_1(R)$ if and only if $I^2=I$. Also $\{0\}$ and $R$ adjacent to all vertices in $\Gamma_1(R)$. We denote the degrees of a vertex $I$ in  $\Gamma_0(R)$ and  $\Gamma_1(R)$,  by $deg_0(I)$ and $deg_1(I)$, respectively. Any loop is counted by multiplicity one in this definition.
\begin{rem}
 \begin{enumerate}
        \item[(i)] The graph $\Gamma_0(R)$ is a null graph i.e., the ring  $R$ has only two ideals if and only if $R$ is a field.
        \item[(ii)] The graph $\Gamma_0(R)$ has only one vertex  i.e., the ring $R$ has only three ideals if and only if $R$ is a local ring with a principal maximal ideal $m=Ra$ such that $a^2=0$. In this case $m$ is also a minimal ideal.
    \end{enumerate}
\end{rem}

Note that the graph $\Gamma_0(R)$ contains comaximal
graph and the complement of intersection graph.
Also if $R$ is a reduced ring, it contains the annihilating ideal
graph. To investigate some properties of $\Gamma_0(R)$, first we state and prove the following lemma:
\begin{lemma}\label{hrd00}
    Let $I$ and $J$ be two ideals of a ring $R$.
    \begin{enumerate}
        \item[(i)] If $I+J=R$, then $IJ=I\cap J$.
        \item[(ii)] If $I\cap J=0$, then $IJ=0$.
        \item[(iii)] If $R$ is a reduced ring, then $IJ=0$ implies $I\cap J=0=IJ$.
        \item[(iv)] If $IJ=J$, then $IJ=I\cap J$.
    \end{enumerate}
\end{lemma}
\proof
            \begin{enumerate}
            \item[(i)] It is obvious that $IJ \subseteq I\cap J$.
            If $I+J=R$, then there are $i\in I$ and $j\in J$ such that $i+j=1$. If  $t\in I\cap J,$ then $t=t(i+j)=ti+tj\in IJ$. So $IJ=I\cap J$.
            \item[(ii)] It follows from  $IJ\subseteq I\cap J$.
            \item[(iii)] It is easy to see that  $(I\cap J)^2\subseteq IJ$. So   $(I\cap J)^2=0$. Since  $R$ is a reduced ring, therefore $I\cap J=0=IJ$.
            \item[(iv)] If $IJ=J$ then $J=IJ\subseteq I$. So $IJ=J=I\cap J$.    \qed

   By Lemma \ref{hrd00} and the definition of $\Gamma_0(R)$ we have the following corollary:

        \end{enumerate}
    \begin{corollary}\label{hrd000}
     \begin{enumerate}
        \item[(i)]  The set of all maximal ideals of $R$,  $Max(R)$ is a clique in $\Gamma_0(R)$.
        \item[(ii)] The set of all minimal ideals of $R$, $M(R)$ is a clique in $\Gamma_0(R)$.
           \end{enumerate}
    \end{corollary}
        We need the following well-known theorem:

    \begin{theorem}\label{hrd04}
        (Nakayama's Lemma) Let $M$ be a finitely generated $R$-module and $I$ be an ideal of $R$. If $IM=M$, then $ann(M)\bigcap (1+I)\neq \emptyset$.
    \end{theorem}

   We recall that a set $S\subseteq V$ is an independent set of a graph $G$, if no two vertices of $S$ are adjacent. The following corollary which is an immediate consequence of Nakayama's Lemma, is useful for determining of independent sets of $\Gamma_0(R)$ (Corollary \ref{hrd004}):
    \begin{corollary}\label{hrd104}
    Let $I$ and $J$ be two ideals of $R$ such that $J$ is a finitely generated ideal and $IJ=J$.
    \begin{enumerate}
        \item[(i)] If  $I\subseteq J(R)$ then $J=0$.
        \item[(ii)] If $ann(J)=0$ i.e.,          $J$ contains a non zero divisor,  then $I=R$.
           \end{enumerate}
    \end{corollary}

    \begin{corollary}\label{hrd004}
     \begin{enumerate}
        \item[(i)] If $\{I_\alpha\}
        $ is a chain of finitely generated proper ideals in $J(R)$, then $\{I_\alpha\}$ is an independent set of $\Gamma_0(R)$.
        \item[(ii)] If $\{I_\alpha\}$ is a chain of finitely generated proper ideals in an integral domain $R$ then $\{I_\alpha\}$ is an independent set of $\Gamma_0(R)$.
           \end{enumerate}
       \end{corollary}
    \proof
        \begin{enumerate}
        \item[(i)] If $I_\alpha\subsetneqq I_\beta$ and $I_\alpha I_\beta=I_\alpha\bigcap I_\beta$ then $I_\alpha=I_\alpha I_\beta$ which is a contradiction by
        part $(i)$ of Corollary \ref{hrd104} (Nakayama's Lemma).
        \item[(ii)] If $I_\alpha\subsetneqq I_\beta$ and $I_\alpha I_\beta=I_\alpha\bigcap I_\beta$ then $I_\alpha=I_\alpha I_\beta$ which is a contradiction by  part $(ii)$ of Corollary \ref{hrd104}.\qed
           \end{enumerate}

    The following lemma is well known and let us  give a proof for it.
    \begin{lemma}\label{hrd0}
        Let $I$ be a finitely generated idempotent ideal of a ring $R$.
        Then $I=Re$ is generated by an idempotent element.
    \end{lemma}
    \proof
    Since $I=I^2$, so $ann(I)\bigcap (1-I)\neq \emptyset$ by Theorem \ref{hrd04}. Hence, there is  $s=1-e\in ann(I)$ such that $sI=0$.
            This implies that $I=Ie$ and $(1-e)e=0$. Therefore $e=e^2$ and
            $I=Re$.\qed

Now we state and prove the following theorem, which is  one of the main result of this section:
     \begin{theorem}\label{hrd1}
        Let $R$ be a ring such that $\Gamma_0(R)$ has order  at least two. The vertex $I$ is adjacent to any other
        vertices if and only if for every  $a\in I$, $a\in Ia$ and
        $I=I^2$.
    \end{theorem}
    \proof
    First assume that vertex $I$ is adjacent to any other
            vertices. We consider two cases:
            \begin{enumerate}
                \item[Case 1)]   $I\bigcap Ann(I)=0$. We have   $0\neq I^2$. If $Ra\subsetneqq I$, then
                $Ra=I\bigcap Ra=Ia\subseteq I^2$. So $a\in Ia$. If $I$ is not a principal
                ideal, then $I=I^2.$ So assume $I=Rb$ is a principal ideal. If
                $0\neq I^2\subsetneqq I$, then $I^3=I^2$. Hence by Lemma
                \ref{hrd0}, we have $I^2=Re$. This implies that $I(1-e)\subseteq Ann(I)\bigcap I=0$.
                Thus $I=Ie\subseteq I^2$ and so  $I=I^2$. Therefore $b\in I^2=Ib$.
                \item[Case 2)]  $I\bigcap Ann(I)\neq 0$. Since $0=IAnn(I)\neq I\bigcap
                Ann(I)$, so $I=Ann(I)$. If $0\neq J\subsetneqq I$, then
                $0=IJ=I\bigcap J=J$ which is a contradiction. So $I$ is a minimal
                ideal. Hence $Ann(I)$ is a maximal ideal and so  $I$ is both a maximal and minimal ideal of $R$. If $R$ has another maximal ideal $m\neq I$, then $0=IAnn(I)=I^2\subseteq m$. So $I\subseteq m$ which is a contradiction. Since every ideal of $R$ sits in a maximal ideal, so  $I$ is the the only nontrivial ideal of $R$ which is a contradiction.
            \end{enumerate}
            Conversely, assume that  for each $a\in I$, $a\in Ia$. If $a\in I\bigcap J$, then $a\in Ia\subseteq IJ$. So $a\in IJ$ and therefore  $IJ=I\bigcap J$.\qed

    \begin{rem}
        Theorem \ref{hrd1} states  that a vertex $I$ is adjacent to all other vertices in $\Gamma_0(R)$ if and only if
        the vertex $I$ is adjacent to all principal ideals contained in I in graph $\Gamma_1(R)$. Also in this case we have a loop in the graph $\Gamma_1(R)$
    \end{rem}
    \begin{theorem}\label{hrd2}
        Let $I$ be a finitely generated ideal  of $R$. Then
        $I$ is adjacent to all other
        vertices if and only if $I=Re$ is generated by an idempotent
    \end{theorem}
    \proof
 Let $J$ be an ideal of $R$ and $t\in J\cap Re$, so $t=re$ and  $t=te\in Je$. Thus  $I$ is adjacent to all other vertices.
            Conversely, If  $I$ is adjacent to all other
            vertices then by Theorem \ref{hrd1}, $I=I^2$.   So  the proof is complete by Lemma \ref{hrd0}.\qed

    \begin{corollary}\label{hrd01}
        If $R$ has a non trivial idempotent, then $\Gamma_0(R)$ is a connected graph and $diam(\Gamma_0(R))\leq 2$.
    \end{corollary}

    The following theorem gives the structure of $\Gamma_0$-graph of Von Neumann regular ring:
    \begin{theorem}\label{hrd6}
    Assume that $\Gamma_0(R)$ has at least two vertices. The ring $R$ is a Von Neumann regular ring if and only if $\Gamma_0(R)(\Gamma_1(R))$ is a complete graph.
    \end{theorem}
    \proof
        If every two vertices are adjacent,  then every principal ideal is
            generated by an idempotent by Theorem \ref{hrd2}. So $\langle
            a\rangle=\langle e\rangle$ where $e$ is an idempotent element and so $a=re$. Thus $ae=re^2=re=a$ and $a=ae^2=eae$. So $R$ is a
            Von Neumann regular ring. Conversely, assume that $R$ is a Von
            Neumann regular ring. If $t\in I\bigcap J$ then, $Rt=Re$ for some
            idempotent $e$. Therefore $t=te\in It\subseteq IJ$ and we have the result.\qed

            The following theorem gives an upper bound for the diameter of $\Gamma_0(R)$ while the Jacobson radical of the ring is zero:
             \begin{theorem}\label{hrd5}
        If $J(R)=0$, then $diam(\Gamma_0(R))\leq 2$.
    \end{theorem}
    \proof
    Let $I,J$ be two distinct ideal of $R$.  If $I\bigcap J=0$ then
            $I$ and $J$ are adjacent by Part $(ii)$ of Lemma \ref{hrd00}.  So
            assume that  $I\bigcap I\neq 0$. Since $J(R)=0$,  there is a maximal
            ideal $m$ such that $I\bigcap J\nsubseteq m$. This implies that
            $I+m=R=J+m$. So $m$ is adjacent to both ideals $I$ and $J$ by
            part $(i)$ of Lemma \ref{hrd00}.\qed

            \medskip

           As an application of Theorem \ref{hrd5}, consider the ring of real continuous functions on a topological space  $X$, i.e., $C(X)$. Since the ideals $M_{x_0}=\{f\in C(X): f(x_0)=0\}$ are  maximal ideals of $C(X)$, so $J(C(X))=0$ and  $diam(\Gamma_0(C(X)))\leq 2$.

        \medskip
        The following theorem is stated in  \cite[p. 15]{shrp} as an exercise.
\begin{theorem}\label{hrd35}
 Let $R$ be  a commutative ring and $f(x)=a_nx^n+\cdots+a_0\in R[x]$.
\begin{enumerate}
\item[(i)] $f(x)$ is a  nilpotent element of $R[x]$ if and only if $a_i$ is a nilpotent element of $R$ for each $i$.
\item[(ii)] $f(x)$ is an invertible element of $R[x]$ if and only if $a_0$ is an invertible element and $a_i$ is nilpotent for each $i\geq 1$.
\item[(iii)] $J(R[x])=Nil(R[x])=Nil(R)[x]$.
\end{enumerate}
\end{theorem}
Now we state and prove the following corollary: 
\begin{corollary}\label{hrd45}
Let $R$ be  a reduced ring. Then $diam(R[x])\leq 2$.
\end{corollary}
\proof 
Since $Nil(R)=0$, so $J(R[x])=Nil(R[x])=Nil(R)[x]=0$. So we have the result by Theorem \ref{hrd5}.\qed

       The following result state a necessary condition for an ideal to be an isolated vertex of $\Gamma_0(R)$.

    \begin{theorem}\label{hrd3}
        Let $R$ be a ring. If $I$ is an isolated vertex of $\Gamma_0(R)$, then $I\subseteq  J(R)$ and $I$ is a large ideal.
    \end{theorem}
    \proof
    Let $m$ be a maximal ideal of $R$. If $I\nsubseteq m$, then
            $I+m=R$. Hence $Im=I\cap m$ which is a contradiction. So
            $I\subseteq J(R).$ If $I\cap J=0$, then $I$ is adjacent to $J$
            which is a contradiction. So $I$ is a large ideal and hence $J(R)$
            is a large ideal.\qed

  The following corollary is an immediate consequence of Theorem \ref{hrd3}:
    \begin{corollary}\label{hrd8}
        If  $\Gamma_0(R)$ is the empty graph, then $R$ is a local ring and every ideal of $R$ is large.
    \end{corollary}

    \begin{theorem}\label{hrd4}
        If $(R,m)$ is  a Noetherian local ring, then  $m$ is an isolated
        vertex of $\Gamma_0(R)$.
    \end{theorem}
    \proof
    If $I\subseteq m$, then by Nakayama's Lemma $Im\subsetneqq
            I=I\cap m$. Hence $m$ is an isolated vertex.\qed

    \begin{corollary}\label{hrd001}
        Let  $R$ be an Artinian ring such that $\Gamma_0(R)$ has at least two vertices. Then  $\Gamma_0(R)$ is a connected graph if and only if $R$ is not a local ring.
    \end{corollary}
    \proof
        First assume that  $\Gamma_0(R)$ is a connected graph. Since every Artinian ring is Noetherian, so the proof is complete by Theorem \ref{hrd4}. Conversely, assume that $R$ is not a local ring. So $R\cong R_1\times\cdots\times R_k$ where $k\geq 2$ by \cite[Theorem 8.7]{shrp}.
        Hence $R$ has a non trivial idempotent. Thus  $\Gamma_0(R)$ is a connected graph by Corollary \ref{hrd01}.\qed

Now, we consider a specific ring and investigate the structure of its $\Gamma_0$-graph:

    \begin{example}
        Let $R=\frac{\mathbb{F}_q[X,Y]}{\langle X,Y\rangle^2}$. It is
        clear that $R$ is an Artinian local ring with $m=\langle
        x,y\rangle$ as maximal ideal. Since $m^2=0$, so every non maximal
        ideal of $R$ correspond to a one dimensional vector subspace of
        two dimensional vector space $m$. So $R$ has $q+1=\frac{q^2-1}{q-1}$ ideal of
        dimension one and a maximal ideal of dimension two. Also one
        dimensional ideals are minimal ideals. Since for every two
        distinct minimal ideals $I$ and $J$ we have $I\bigcap J=0$, so
        every two distinct minimal ideals $I$ and $J$ are adjacent in
        $\Gamma_0(R)$. So $\Gamma_0(R)$ is union of a complete graph
        $K_{q+1}$ and $K_1$, i.e., $\Gamma_0\big(\frac{\mathbb{F}_q[X,Y]}{\langle X,Y\rangle^2}\big)=K_{q+1}\cup K_1$.
    \end{example}

The following theorem describes the $\Gamma_0$-graph for discrete valuation ring:

    \begin{theorem}\label{hrd7}
        If $R$ is a discrete valuation ring(DVR), then $\Gamma_0(R)$ is the empty graph.
    \end{theorem}
    \proof
Suppose that  $R$ is a discrete valuation ring(DVR). It is well known that $R$ is a local ring such that its maximal ideal $m$ is principal and only ideals of $R$
            are $m^i$. If $m^i$ and $m^j(j>i)$ are adjacent then $m^im^j=m^i\cap m^j=m^j$. So Nakayama's Lemma implies that
            $m^j=0$ which is a contradiction. Hence  $\Gamma_0(R)$ is the empty graph.\qed

 To study the structure of $\Gamma_0(\mathbb{Z}_n)$, we need the following result.

            \begin{theorem}\label{hrd010}
            	Let  $R\cong R_1\times\cdots\times R_k$. Then  $\Gamma_1(R)\cong \Gamma_1(R_1)\times\cdots\times \Gamma_1(R_k).$
            \end{theorem}
            \proof
            It is well known that every ideal $I$ of $R$ is equal to $I_1\times\cdots\times I_k$ for some ideals $I_i$ of $R_i$.
            Also the vertices  $I=I_1\times\cdots\times I_k$ and $J=J_1\times\cdots\times J_k$ are adjacent if and
            only if $I_iJ_i=I_i\bigcap J_i$ for each $i$. So the proof is complete.\qed
            \begin{rem}\label{hrd11}
            	If $R=\mathbb{Z}_{p_1^{\gamma_1}}\times\cdots\times\mathbb{Z}_{p_k^{\gamma_k}}$($p_i$'s are not necessarily distinct), then ideals of $R$ are
            	$\langle(p_1^{\alpha_1},\ldots, p_k^{\alpha_k})\rangle$ where $0\leq \alpha_i\leq \gamma_i$. Also two ideals
            	$I=\langle(p_1^{\alpha_1},\ldots, p_k^{\alpha_k})\rangle$ and $J=\langle(p_1^{\beta_1},\ldots, p_k^{\beta_k})\rangle$ are adjacent if and only if $max\{\alpha_i,\beta_i\}=min\{\alpha_i+\beta_i,\gamma_i\}$. So we can construct
            	$\Gamma_1(R)$ as follows:
            	
            	The vertex set is $\{(\alpha_1,\ldots,\alpha_k):0\leq \alpha_i\leq \gamma_i\}$ and two vertices
            	$(\alpha_1,\ldots,\alpha_k)$ and $(\beta_1,\ldots,\beta_k)$ are adjacent if and only if $max\{\alpha_i,\beta_i\}=min\{\alpha_i+\beta_i,\gamma_i\}$.
            \end{rem}
            	\begin{example}
            		Consider $R=\mathbb{Z}_{20}$. We have $\mathbb{Z}_{20}=\mathbb{Z}_{2^2}\times \mathbb{Z}_5$. We shall draw $\Gamma_0(\mathbb{Z}_{20})$. The vertex set of this graph is
            		$V(\Gamma_0(\mathbb{Z}_{20}))=\{(1,0),(2,0),(0,1),(1,1)\}$. By Remark \ref{hrd11} we have the Figure \ref{figure1} for this graph:
            		
            		\begin{figure}[!ht]
            			\begin{center}
            				\includegraphics[]{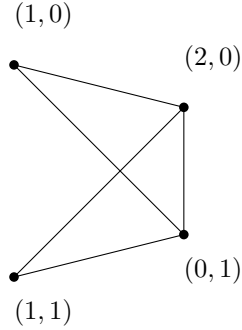}
            				\caption{ \label{figure1} The graph
            					$\Gamma_0(\mathbb{Z}_{20})$.}
            			\end{center}
            		\end{figure}
            	\end{example}

            	The following theorem gives the degree of the  vertices of $\Gamma_0(\mathbb{Z}_n)$.
            
            \begin{theorem}\label{hrd12}
            	Let $n=p_1^{\gamma_1}\cdots p_k^{\gamma_k}$ and $a=p_1^{\alpha_1}\ldots p_k^{\alpha_k}$ where $0\leq \alpha_i\leq \gamma_i$ be a divisor of $n$. Set $R=\mathbb{Z}_n$ and  $I=\langle a\rangle$. Suppose that
            	$A(I)=\{i:1\leq i\leq k,0\leq\alpha_i\leq\gamma_i\}$. The degree of ideal $I$ in $\Gamma_0(R)$ is
            	$$deg_0(I)=\Big(2^{|A(I)|}\prod_{i\notin A(I)}(\gamma_i+1)\Big)-2-\Big\lfloor\frac{1}{|A(I)|+1}\Big\rfloor.$$
            \end{theorem}
            \proof
            According to Remark \ref{hrd11}, we do computations with $(\alpha_1,\ldots,\alpha_k)$. Let $B=N_{\Gamma_1}(I)$ be the open neighborhood of the vertex $I$. So $N_{\Gamma_0}(I)=B\setminus\{1,I,n\}$.
            Assume $b=(\beta_1,\cdots,\beta_k)\in B$. Hence
            $max\{\alpha_i,\beta_i\}=min\{\alpha_i+\beta_i,\gamma_i\}$. This
            implies that if $i\notin A(I)$, then
            we have the result. If $i\in A(I)$ then $\beta_i=0$ or $\beta_i=\gamma_i$.
            So in the first case $\beta_i$ can be any number of the set $\{0,\cdots,\gamma_i\}$. In the last case $\beta_i=0$ or $\beta_i=\gamma_i$.
            So we can  choose  $b$  in  $2^{|A(I)|}\prod_{i\notin A(I)}(\gamma_i+1)$ ways. Two of these $b$ correspond to $(0,\ldots,0)$
            and $(\gamma_1,\ldots,\gamma_k)$.
            If $I\in B$ then $min\{2\alpha_i,\gamma_i\}=\alpha_i$. So $\alpha_i=0$ or $\alpha_i=\gamma_i$. Thus $A(I)=\emptyset$. Conversely,
            if $A(I)=\emptyset$ then for each $i,\; min\{2\alpha_i,\gamma_i\}=\alpha_i$. Hence $I\in B$.
            So we must exclude $I$ from $B$ in this case. So $deg_0(I)=\Big(2^{|A(I)|}\prod_{i\notin A(I)}(\gamma_i+1)\Big)-2-\lfloor\frac{1}{|A(I)+1}\rfloor$.\qed
            \begin{example}
            	Let $n=36=2^23^2$. So $\mathbb{Z}_{36}\cong \mathbb{Z}_4\times\mathbb{Z}_9$ and $\gamma_1=\gamma_2=2$.
            	If $I=6=2\times 3$ then $A(6)=\{1,2\}$ and $|A(6)|=2$. So
            	$deg_0(6)=2^2-2=2$. The neighbors of the vertex $6$ are $4$ and $9$. If $I=4$ then $A(4)=\emptyset$ and $|A(4)|=0$. So
            	$deg_0(4)=3\times 3-2-1=6$. Note that $\langle 4\rangle=\langle 16\rangle$ i.e. $\langle 4\rangle$ is an idempotent ideal. So $N_{\Gamma_0}(4)=\{2,3,6,9,12,18\}$. We have shown
            	$\Gamma_0(\mathbb{Z}_{36})$ in Figure \ref{2}.

            	\begin{figure}[!ht]
            		\begin{center}
            			\includegraphics[]{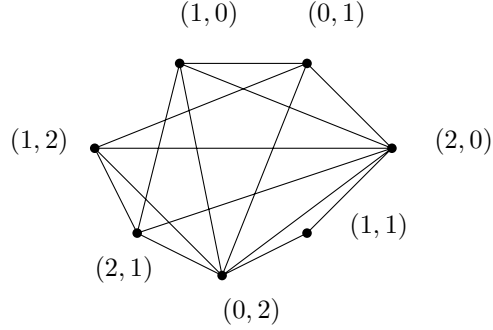}
            			\caption{ \label{2} The graph $\Gamma_0(\mathbb{Z}_{36})$.}
            		\end{center}
            	\end{figure}
            	
            \end{example}

\section{Primary decomposition of ideals of $R$ and connectivity of $\Gamma_0(R)$}

In this section we shall study the relationship between the primary decomposition of ideals of a ring  $R$ and the  connectivity of $\Gamma_0(R)$. We begin with the following  result.

 \begin{lemma}\label{hrd0'}
 	(Prime avoidance lemma\cite{shrp})
 	Suppose that
 	$I\subseteq\cup_{i=1}^n P_i$, where $P_i$'s are prime ideals. Then  $I\subseteq P_i$ for some $1\leq i\leq n$.
 \end{lemma}
 \begin{theorem}\label{hrd6'}
 	Let $I$ be an ideal of the ring $R$. Then $I\bigcap Ra=Ia$ if and only if $(I:a)=I+ann(a)$. In particular, if $(I:a)=I$ then $I\bigcap Ra=Ia$.
 \end{theorem}
 \proof
 Assume that  $x=ra\in Ra\bigcap I$. So $r\in (I:a)$. Hence $r=i+s$ where $i\in I$ and $s\in ann(a)$.
 		Thus $x=ra=ia+sa=ia\in Ia$. Conversely, assume that $r\in (I:a)$. So $x=ra\in I\bigcap Ra=Ia$.
 		hence there is an $i\in I$ such that $ra=ia$. This implies that $r-i\in ann(a)$. So $r=i+(r-i)\in I+ann(a)$.\qed

 \begin{theorem} \label{hrd05}
 	\begin{enumerate}
 		\item[(i)] 	Let $Q$ be  a primary ideal and $a\notin r(Q)=P$. Then $(Q:a)=Q$.
 		
 		\item[(ii)]  Suppose that  $I$ has a primary decomposition.  If $a\notin \bigcup_{P\in Ass(I)} P$, then $(I:a)=I$. 		
 		 	\end{enumerate}
 \end{theorem}
 \proof
 \begin{enumerate}
 \item[(i)]
 	Let $b\in (Q:a)$. So $ba\in Q$. If $b\notin Q$ then $a\in r(Q)$ which is a contradiction.
 	\item[(ii)] Let $I=\bigcap_{i=1}^n Q_i$. Then $(I:a)=(\bigcap_{i=1}^n Q_i:a)=\bigcap_{i=1}^n (Q_i:a)=\bigcap_{i=1}^n Q_i=I$.\qed
 	
 	 \end{enumerate}
 	
 	Since every prime ideal is a primary ideal, we have the following corollary:

 \begin{corollary}\label{hrd1'}
 	Let $P$ be a prime ideal of a ring $R$. If  $a\notin P$,  then $Ra\bigcap
 	P=Pa=PRa$, i.e., $P$ and $Ra$ are adjacent.
 \end{corollary}

 \begin{corollary}\label{hrd2'}
 	Suppose that $R$ is not an integral domain. If $ab\notin Nil(R)$ then
 the distance between $Ra$ and $Rb$ is not more than two, i.e,	$d(Ra,Rb)\leq 2$.
 \end{corollary}
 \proof
 Since $Nil(R)$ is the intersection of all prime ideals, so there
 		is a prime ideal $P$ such that $ab\notin P$. Hence $a,b\notin P$.
 		So $Ra$ and $Rb$ are adjacent to $P$ by Corollary \ref{hrd1'}.\qed

 \begin{theorem}\label{hrd30}
 	Suppose that two ideals  $I$ and $J$  have primary decomposition. If $Max(R)\nsubseteq
 	(Ass(I)\bigcup Ass(J))$, then $d(I,J)\leq 2$.
 \end{theorem}
 \proof
 Let $m\in Max(R)\backslash
 		(Ass(I)\bigcup Ass(J))$. So there is an $a\in m\backslash \bigcup_{P\in Ass(I)\bigcup Ass(J)} P$ by prime avoidance lemma. So $I$ and $J$  are adjacent
 		to $Ra$ by Theorem \ref{hrd05}(ii) and the proof is complete.\qed
 		
 \begin{corollary}\label{hrd3'}
 	\begin{enumerate}
 	\item[(i)]
 			Let $P,Q\in Spec(R)$ be  two prime ideals. If $Max(R)\nsubseteq
 		\{P,Q\}$ then $d(P,Q)\leq 2$.
 		
 		\item[(ii)] 	Suppose that  $R$  has at least three maximal ideal. Then for every two
 		prime ideals $P$ and $Q$,  $d(P,Q)\leq 2$.
 	\end{enumerate}

 \end{corollary}

 \begin{theorem}\label{7}
 	Let $R$ be a Noetherian ring with infinitely many maximal ideals. Then  $diam(\Gamma_0(R))\leq 2$.
 \end{theorem}
 \proof
 Let $I$ and $J$ be two ideals of $R$. It is well known that every ideal in a Noetherian ring has a primary decomposition (\cite{shrp}). So the proof is complete by Theorem \ref{hrd30}.\qed

\medskip
 The following example shows that the condition $|Max(R)|=\infty$ is necessary in the Theorem \ref{7}.
 \begin{example}
 	Let $p_1,\ldots,p_k$ be distinct prime numbers. Set $S=\mathbb{Z}\backslash\bigcup_{i=1}^k p_i\mathbb{Z}$ and $R=S^{-1}\mathbb{Z}$. Then $R$ is a PID with $k$ distinct prime $p_1,\ldots,p_k$. If $I=\langle a\rangle$
 	and $J=\langle b\rangle$ are adjacent then $\langle a,b\rangle=R$. So if $I\subseteq J(R)=p_1\cdots p_kR$ then $I$ is an isolated vertex.
 \end{example}

\medskip
 Here, we state and prove the following lemma to obtain the diameter of $\Gamma_0$-graph of the polynomial ring:

 \begin{lemma}\label{polyn}
 	Let $R$ be a commutative ring. The polynomial ring $R[x]$ has infinitely many maximal ideals.
  \end{lemma}
  \proof
  Let $M$ be a maximal ideal of $R$ and $F=\frac{R}{M}$ be its residue field. Since $F[x]\cong \frac{R[x]}{M[x]}$, it suffices to prove the lemma for $F[x]$. Every maximal ideal of $F[x]$ is generated by an irreducible polynomial, because, $F[x]$ is a PID. If $F$ is an infinite field, the set $\{\langle x-a\rangle:a\in F\}$ is an infinite set of maximal ideals.  Now suppose that $F$ is a finite field, then by a well-known result, for every $n\in \mathbb{N}$, there is an irreducible polynomial of degree $n$. Therefore we have the result.\qed

  By Theorem \ref{7} and Lemma \ref{polyn}, we have the following corollary:

  \begin{corollary}
  If $R$ is a Noetherian ring, then $diam(\Gamma_0(R[x]))\leq 2$.
  \end{corollary}

 By  checking the $\Gamma_0$-graph for  well-known rings, we think that the diameter  of every connected component of $\Gamma_0(R)$ is not more than two. So, we end the paper by the following conjecture:
\begin{conjecture}
The diameter of every connected component of $\Gamma_0(R)$ is not more than two.
\end{conjecture}

\end{document}